\documentclass[a4,12pt,reqno,fleqn]{amsart}
\usepackage{graphicx}
\usepackage{amscd}
\usepackage{xfrac}
\usepackage{mathtools}

\emergencystretch=\maxdimen
\makeatletter
\@namedef{subjclassname@2020}{\textup{2020} Mathematics Subject Classification}
\makeatother



\usepackage{amsmath}
\usepackage{amsthm}
\usepackage{thmtools}
\usepackage{amsfonts,amssymb,hyperref,mathrsfs,cleveref,setspace,enumitem,verbatim}
\usepackage[margin=3.5cm]{geometry}

\allowdisplaybreaks[1]
\usepackage{array}
\usepackage{xcolor}
\usepackage[all]{xy}

\DeclareMathAlphabet{\mathpzc}{OT1}{pzc}{m}{it}

\newtheorem{thm}{Theorem}[section]
\newtheorem{cor}[thm]{Corollary}
\newtheorem{lem}[thm]{Lemma}
\newtheorem{propn}[thm]{Proposition}
\theoremstyle{definition}
\newtheorem{defn}[thm]{Definition}
\newtheorem{remark}[thm]{Remark}
\newtheorem{example}[thm]{Example}

\newcommand{\seq}{\subseteq}
\newcommand{\oh}{\otimes_h}
\newcommand{\omin}{\otimes_{\min}}
\newcommand{\op}{\otimes_{\gamma}}
\newcommand{\oa}{\otimes_{\alpha}}

\newcommand{\ra}{\rightarrow}
\newcommand{\ot}{\otimes}

\newcommand{\C}{\mathbb{C}}

\newcommand{\mcal}{\mathcal}
\author[A. Paliwal and R. Jain]{Anmol Paliwal and Ranjana Jain}
\address{Anmol Paliwal, Department of  Mathematics, University of Delhi, Delhi}
\email{anmol.paliwal1234@gmail.com}

\address{Ranjana Jain, Department of Mathematics, University of Delhi, Delhi}
\email{rjain@maths.du.ac.in}
\thanks{Research of the first named author is supported by NBHM Fellowship vide S.No. 0203/7/2023/R\&D-11/6429. \\
 Conflict of Interest Statement: The authors declare no conflicts of interest regarding this manuscript. \\
 Data availability:  No new data were created or analysed in this study. Data sharing is not applicable to this article.}
\begin{document}
	
	\title[Weak Centrality for Certain Tensor Products of $C^\ast$-algebras]{Weak Centrality for Certain Tensor Products of $C^\ast$-algebras}
	\maketitle
	
	\begin{abstract}
	In this article, we discuss the weak centrality of the tensor product $A\otimes_\alpha B$ of $C^\ast$-algebras $A$ and $B$ in terms of the weak centrality  of $A$ and $B$, where $\alpha$ is either the Haagerup  or the Banach space projective tensor product. In  the due course, we also identify the largest weakly central ideal of $A\otimes_\alpha B$ in certain cases. Centralilty and quasi-centrality of these tensor products are also discussed.  
	
	\end{abstract}

	\section{Introduction}
	
	 Kaplansky \cite{Kaplansky}, in 1948, introduced the notion of central $C^*$-algebras in order to study those $C^*$-algebras to which commutative methods can be generalized. A $C^\ast$-algebra $A$ is called \textit{central} if two primitive ideals $P$ and $P'$ of $A$ coincide whenever $P\cap Z(A)=P'\cap Z(A)$, 
	 where $Z(A)$ is the centre of $A$. Later, in 1951, in order to study an analogue of Dixmier trace, this concept was modified by Misonou and Nakamura \cite{MisonouNakamura} by introducing  weakly central $C^*$-algebras, where the primitive ideals in the definition of a central algebra were replaced by maximal ideals. 	 Misonou \cite{Misonou} further proved that any unital weakly central $C^\ast$-algebra can be decomposed into factorial $C^\ast$-algebras, and that every von Neumann algebra is weakly central. 
	 
	 In the last few decades,  operator algebraists have done remarkable work on weakly central $C^*$-algebras (see \cite{Archbold,ArchboldRobertTikuisis,ArchboldGogic}).  Archbold et al. \cite{ArchboldRobertTikuisis} proved that a $C^\ast$-algebra $A$ satisfying Dixmier property is weakly central, 
	 whereas the converse is true  for a postliminal weakly central $C^\ast$-algebra. Archbold \cite[Theorem 3.1]{Archbold} also proved that two unital $C^*$-algebras $A$ and $B$ are weakly central if and only if $A\otimes_\beta B$ is weakly central for some (or, every) $C^\ast$-tensor norm $\beta$. Archbold and Gogic  also discussed the largest weakly central ideal of a $C^*$-algebra in detail.  In particular, they proved the existence, and provided a description of the  largest weakly central ideal in any $C^\ast$-algebra $A$. Interestingly, the largest weakly central ideal turns out to be $\ker T_A :=\cap_{M \in T_A} M$,  where $T_A$ is the collection of those maximal modular ideals which prevent $A$ from being weakly central.
	 
	 On the other hand, the theory of tensor products of $C^*$-algebras  is an  indispensable tool in the proper  understanding of these objects.  It is therefore quite natural to ask whether the various notions of centrality of individual algebras are related to those of their tensor products. In this article, we shall ponder upon these questions for the Haagerup tensor product ($\oh$) and the Banach space projective tensor product ($\op$) of two $C^*$-algebras. The fact that $A \oh B$ and $A \op B$ are Banach and Banach $*$-algebras, respectively, and these are $C^*$-algebras only in trivial cases makes the theory more intriguing. We first introduce the notion of weakly central (respectively, central) Banach algebras, and study their properties and some examples.	We prove that for $C$*-algebras $A$ and $B$, $A\otimes_\alpha B$ is weakly central if $A$ and $B$ are weakly central, and the converse is true if $Max(A)$ and $Max(B)$ both are non-empty, where $\oa = \oh$ or $\op$, and $Max(A)$ represents the collection of maximal modular ideals of $A$. We then move on to discuss whether we can express the largest weakly central ideal of $A \oa B$ in terms of those of $A$ and $B$. To our satisfaction, we obtain affirmative answer in certain cases by describing $\ker T_{A\oa B}$ in terms of $\ker T_A$ and $\ker T_B$.  
	  Lastly, we investigate centrality and quasi-centrality for these tensor products. 
	
		
	\section{Some properties of the Banach space projective tensor product of $C^\ast$-algebras}
In this section, we shall collect certain properties related to the Banach space projective tensor product of $C^*$-algebras which will be needed for the main results. We would like to mention that some of these results are the analogues of Haagerup product (see \cite{AllenSinclairSmith, ArchboldKaniuthSchlitingSomerset}) and operator space projective tensor product of $C^*$-algebras (see \cite{JainKumar}),  and the proof of these results follow on the similar lines. However, we present the proofs of these results for the sake of completeness and for the convenience of the reader. 

	Let us fix some notations which will be useful throughout the article.
If $X$ is a Banach algebra, $Id(X)$ represents the collection of closed ideal of $X$. For $S\subseteq Id(X), \; \ker S$ represents the intersection of members of $S$. An ideal $I$ of a Banach algebra $X$ is said to be {\it modular} if $X/I$ is unital, or equivalently, if there exists an $u\in X$ such that $iu-u, \ ui-u\in I$ for all $i\in I$.
We denote $\mathcal{M}(X)=\ker Max(X)$. 

For $C^*$-algebras $A$ and $B$, the {\it Banach space projective tensor norm} and the {\it Haagerup tensor norm} of $ u \in A\otimes B$ are defined as
$$||u||_{\gamma}=\inf \left\{\sum\limits_{i=1}^n{||a_{i}||||b_{i}||}: u=\sum\limits_{i=1}^n{a_{i}\otimes b_{i}}\right\},  \text{ and}$$
 $$||u||_h =\inf \left\{ \bigg\| \sum\limits_{i=1}^n a_{i}a_{i}^{\ast}\bigg\|^{1/2}\bigg\|\sum\limits_{i=1}^{n} b_{i}^{\ast}b_{i}\bigg\|^{1/2}: u=\sum\limits_{i=1}^n{a_{i}\otimes b_{i}} \right\}.$$
The completion of the algebriac tensor product $A \otimes B$ under these norms are denoted by $A\op B$ and $A \oh B$, respectively, and are known as the {\it Banach space projective tensor product} and the {\it Haagerup tensor product} of $A$ and $B$, respectively. It is well known that $A \oh B$ and $A\op B$ are Banach and  Banach*-algebra, respectively.  Throughout the article, $A$ and $B$ will denote $C^\ast$-algebras. 

	For $\varphi \in A^*$, the right slice map $R_{\varphi}:A
	\otimes B \rightarrow B$ is given by
	$$ R_{\varphi}\left(\sum_1^n a_i \otimes b_i\right)=\sum_{1}^{n}
	\varphi (a_i) b_i.$$ Since $ \|R_{\varphi} ( \sum_1^n a_i \otimes b_i
	) \|  \leq \| \varphi\| \sum_1^n  \| a_i \| \| b_i\| \leq \| \varphi\| \|\sum_1^n a_i \otimes b_i\|_{\gamma}$, it
	follows that $R_\varphi$ is continuous with respect to  $ \| \cdot
	\|_\gamma$ with $\|R_{\varphi} \| \leq \|\varphi\|$.  So we can
	extend this map to a bounded linear map $R_{\varphi}: A \op B
	\rightarrow B$. Similarly, the left slice map $L_{\psi} \in B(A\op
	B, A)$ for every $\psi \in B^*$. For each closed ideal $J$ in $B$, $A\op J$ is a closed ideal in
	$A\op B$ and clearly $R_{\varphi}(A \op J) \seq J$ for all $\varphi \in
	A^*$. We next prove the converse of this statement, which is also known as the slice
	map problem (for ideals) (see \Cref{slice}).

		Note that for $\varphi \in A^*, \, \psi \in B^*$, and for  $\sum{a_i \otimes b_i} \in A
	\otimes B$, we have
	\begin{eqnarray*}
		\bigg|(\varphi \otimes \psi)\left(\sum{a_i \otimes b_i}\right)\bigg| & = &
		\bigg|\sum{\varphi(a_i) \psi(b_i)}\bigg|\\
		&\leq & \| \varphi \| \| \psi \| \bigg\|\sum{a_i \otimes b_i} \bigg\|_{\lambda}\\
		&\leq & \| \varphi \| \| \psi \| \bigg\|\sum{a_i \otimes b_i} \bigg\|_{\min}\\
		 & \leq & \| \varphi \| \| \psi \| \bigg\|\sum{a_i \otimes b_i} \bigg\|_{\gamma},		 
	\end{eqnarray*}
	where $\| \cdot \|_{\lambda}$ is the Banach space injective tensor norm and $\| \cdot \|_{\min}$ is the injective tensor norm of $C^\ast$-algebras. This implies that 
		  $\varphi \otimes \psi$ can be extended
	continuously to a functional $
	\varphi \omin \psi$ on $ A\omin B$ and to $\varphi \op \psi$ on $A\op B$. Recall that for a Banach space $E$, a subset $S$ of $E^*$ is said to
	be {\it total} on $E$ if for $x \in E$, $\theta (x) = 0$ for all $\theta
	\in S$ implies that $x= 0$.
	
	\begin{propn}\label{mintotal}
		The set $S:= \{ \varphi \omin \psi : \varphi \in A^*, \, \psi \in
		B^* \}$ is total on $A \omin B$. 
	\end{propn}
	\begin{proof}
		Let $\{\pi_A,H\}$ and $\{\pi_B,K\}$ be faithful representations of $A$ and $B$, respectively. Set $\varphi(a)=\langle \pi_A(a)h_1,h_2\rangle$ and $\psi(b)=\langle \pi_B(b)k_1,k_2\rangle$ for $h_i\in H$ and $k_i\in K$. Then $\varphi\in A^\ast$ and $\psi\in B^\ast$ and $(\varphi \omin \psi)(x)=\langle \pi_A\otimes \pi_B(x)h_1\otimes k_1,h_2\otimes k_2\rangle$. If  $(\varphi \omin \psi)(x)=0$ for all $\varphi \in A^*$ and $\psi \in B^*$, then $x=0$ since $\pi_A\otimes \pi_B$ is faithful \cite[Theorem IV.4.9]{Takesaki}. Hence, the set $S$ is total.
	\end{proof}
	Using this, we obtain the following.
	
	\begin{cor}\label{obp-total}
		The set $S=\{\varphi \op \psi : \varphi \in A^*, \, \psi \in B^* \}$
		is total on $A \op B$.
	\end{cor}
	
	\begin{proof}
		 Let $x \in A\op B$ be such that $(\varphi \op
		\psi )(x)=0$ for all $\varphi \in A^*$ and $\psi \in B^*.$ Choose a
		sequence $\{x_n\}$ in $A\otimes B$ such that $\|x_n - x\|_{\gamma}
		\rightarrow 0$. Then, $$
		\lim_n (\varphi \otimes \psi)(x_n)=\lim_n (\varphi
		\op \psi)(x_n)= (\varphi \op \psi)(x)=0.$$
		It is well known that the identity map on $A \ot B$ extends
		to an injective $*$-homomorphism $i: A\op B \rightarrow
		A\omin B$ (see \cite{haag}). Since $\| \cdot\|_{\min} \leq \| \cdot \|_{\gamma}$, $i$ is a
		contraction and we see that $\| x_n - i(x) \|_{\min} \ra 0$ as
		well. This implies that
		$$(\varphi \omin \psi)(i(x))=\lim_n (\varphi \omin
		\psi)(x_n)=\lim_n (\varphi \otimes \psi)(x_n)= 0,$$ 
		and this true for
		all $\varphi \in A^*$ and $\psi \in B^*$. Therefore, by
		\Cref{mintotal}, $i(x)=0$, and by injectivity of $i$, we have $x=0$ as well.
	\end{proof}

	\begin{cor}\label{total} 
		The set \{$R_\varphi : \varphi \in A^*$\} is total on $A\op B$.
	\end{cor}
	
	\begin{proof}
	Note that for $ x \in A \op B$,		$$
			\langle x, \varphi \op \psi \rangle = \langle R_{\varphi}(x),
			\psi\rangle = \langle L_\psi (x), \varphi\rangle \,
			$$
			 where $\langle
		\cdot, \cdot \rangle$ denotes the duality bracket. The result now follows from	\Cref{obp-total}.
	\end{proof}

	\begin{propn}\label{slice}
		Let $J$ be a closed ideal in $B$ and
		$x\in A \op B$ be such that $R_{\varphi}(x) \in J$ for all $\varphi \in
		A^*$. Then $x \in A\op J$. In other words,
		$$A\op J= \{ x \in A\op B : R_\varphi (x) \in J \ \mathrm{for\ all}\ \varphi \in A^*\}.$$
	\end{propn}
	
	\begin{proof}
		Let $K= \{ x \in A\op B : R_\varphi (x) \in J
		\ \mathrm{for\ all}\ \varphi \in A^*\}$ and let $x\in K$. By
		\cite[Proposition 2.5]{ryan}, the tensor map $\text{Id} \ot \pi$
		extends to a quotient map $ \text{Id} \op \pi: A\op B \rightarrow
		A\op B/J$, where $\text{Id}$ is the identity map on $A$ and $\pi: B
		\ra B/J$ is the canonical quotient map. Then it is easily
		seen that $ \pi \circ R_\varphi = r_{\varphi} \circ (i \op \pi)$
		for all $\varphi \in A^*$, where $r_\varphi:A\op B/J \rightarrow
		B/J$ is the right slice map.  Since $R_\varphi (x) \in J$ for all
		$\varphi \in A^*$, we see that $\pi(R_\varphi (x)) = 0$ for all
		$\varphi \in A^*$, and then the above relation gives $r_\varphi (
		(\text{Id} \op \pi) (x)) = 0$ for all $\varphi \in A^*$. Therefore,
		by \Cref{total}, we have $(\text{Id} \op \pi) (x) = 0$. Thus, by \cite[Lemma 1]{GuptaJain Banach}, 
		$x \in \ker(\text{Id} \op \pi) = A \op J$, proving the claim.
	\end{proof}

	The above results are also true for the left slice map, with their
	corresponding statements. In particular, we have:
	
	\begin{propn}\label{slice1}
		Let $I$ be a closed ideal in $A$ and
		$x\in A\op B$ be such that $L_{\psi} (x) \in I$ for all $\psi
		\in B^*$. Then $x \in I\op B$.
	\end{propn}
	
	We can now prove that the intersection of product ideals is again a product ideal.
	\begin{thm}\label{intersection}
		Let $J_i$ and $K_i$, $i=1,2$ be closed ideals of $A$ and $B$, respectively. Then,
		$$ (J_1 \op K_1) \cap (J_2 \op K_2) = (J_1 \cap J_2) \op (K_1 \cap K_2).$$
	\end{thm}
	
	\begin{proof}
		First of all, it is easy to see that $$ (J_1 \cap J_2) \op
		(K_1 \cap K_2) \subseteq (J_1 \op K_1) \cap (J_2
		\op K_2). $$ For the reverse containment, consider an element
		$v \in (J_1 \op K_1) \cap (J_2 \op K_2)$ and
		$\varphi \in A^*$. Now $v\in J_1 \op K_1$, so $R_\varphi(v) \in
		K_1$ (restriction of $\varphi$ on $J_1$). Similarly $v\in J_2
		\op K_2$, so $R_\varphi(v) \in K_2$, which implies $R_\varphi(v)
		\in (K_1 \cap K_2)$ and this is true for any $\varphi \in A^*$, so by
		\Cref{slice},  $v\in A \op (K_1 \cap K_2)$.
		
		Next, consider any $\psi \in (K_1 \cap K_2)^*$ and consider its
		extension $\widetilde{\psi}$ on $B^*$. Then again using the above argument
		with left slice map, we get $L_{\widetilde{\psi}}(v) \in (J_1 \cap
		J_2)$. But note that $L_\psi(v) = L_{\widetilde{\psi}}(v)$, so that
		$L_\psi(v) \in (J_1 \cap J_2)$. This is true for every $\psi \in
		(K_1 \cap K_2)^*$ and $v\in A\op (K_1 \cap K_2)$. Combining all these
		facts together with \Cref{slice1}, we obtain that $v\in (J_1 \cap J_2)
		\op (K_1 \cap K_2)$, proving the claim.
	\end{proof} 
	
	\begin{cor}
		Let $I$ and $J$ be closed ideals in $A$ and $B$, respectively. Then, 
		$$ I \op J = \{ x \in A \op B: R_{\varphi} (x ) \in J\,
		\mathrm{and}\ L_{\psi} (x) \in I \ \mathrm{for\ all}\ \varphi \in A^*,
		\psi \in B^*\}.$$ \end{cor}
	
	Let $\mcal{C}$ denote the collection of finite sums of product ideals
	in $A \op B$.  We next show that finite intersection of
	finite sums of product ideals is again a finite sum of product ideals.
	
	\begin{cor}
		$\mcal{C}$ is closed under finite intersection.
	\end{cor}
	\begin{proof}
		By \cite[Proposition 4.11]{GuptaJain Nach}, every $I \in \mcal{C}$ contains a
		quasi-central approximate identity. So, if $I = \sum_{r = 1}^k I_r
		\op J_r, I' = \sum_{s = 1}^l I_s' \op J_s' \in \mcal{C}$, then by
		\cite[Proposition 3.6]{AllenSinclairSmith}, $I \cap I' = \sum_{r, s} (I_r \op J_r) \cap
		(I_s' \op J_s')$. Therefore, by \Cref{intersection}, we have $$I \cap I' =
		\sum_{r, s} (I_r \cap I_s') \op (J_r \cap J_s'),$$
		and hence the result.
		\end{proof}
	\begin{propn}\label{kerideals}
		For $U\subseteq Id(A)$ and $V\subseteq Id(B)$,
		$$A\op \big(\bigcap_{N\in V} N\big) + \big(\bigcap_{M\in U} M\big) \op B=\bigcap\{A\op N+M\op B| \; M\in U, N\in V\}.$$
	\end{propn}	
	\begin{proof}
		Let us denote $L=\bigcap \{A\op N+M\op B| \; M\in U, N\in V\}$, $M_0=\ker(U)$ and $N_0=\ker(V)$, then it is trivial that $A\op N_0+M_0\op B\subseteq L.$
		Consider the map $q_{M_0}\otimes q_{N_0} : A\op B\rightarrow (A/M_0)\op (B/N_0)$ where $q_{M_0}$ and $q_{N_0}$ are canonical quotient maps. If the containment is strict then $(q_{M_0}\otimes q_{N_0})(L)$ is a non-zero closed ideal in $(A/M_0)\op (B/N_0)$. By \cite[Corollary 2]{GuptaJain Banach}, there exist $a\in A\backslash M_0$ and $b\in B\backslash N_0$ such that $a\otimes b\in L$. Since $a\not\in M_0, b\not\in N_0$ there exist $M\in U$ and $N\in V$ such that $a\not\in M, b\not\in N$. Let $\phi\in A^*$ and $\psi\in B^*$ such that $\phi(a)\neq 0, \psi(b)\neq 0$ but $\phi(M)=0=\psi(N)$. It implies $a\otimes b\not\in A\op N+M\op B$ as $(\phi \otimes \psi) (a\otimes b)\neq 0$ but $(\phi \otimes \psi) (A\op N+M\op B)=0$. Thus, $a\otimes b \not\in L$ which is a contradiction.
	\end{proof}
	 Next, we characterise the primitive ideals of the projective tensor product. Here, $Prim(A)$ and $Prime(A)$ represent the collection of primitive ideals and prime ideals of $A$, respectively.
	 
	\begin{lem}\label{Factorrepresent}
		Let $\pi$ be an irreducible $*$-representation of $A\otimes_\gamma B$ on a Hilbert space $H$. Then, there exist factor $*$-representations $\pi_1$ and $\pi_2$ of $A$ and $B$ on $H$ with commuting ranges such that 
		$$\pi(a\otimes b)=\pi_1(a) \pi_2(b)\; for \; all \; a\in A, \; b\in B.$$
		\end{lem}
		\begin{proof}
			By \cite[Lemma IV.4.1]{Takesaki},  there exist *-representations $\pi_1$ and $\pi_2$ of $A$ and $B$ on $H$ with commuting ranges such that 
			$$\pi(a\otimes b)=\pi_1(a)\pi_2(b)\; \text{ for all} \; a\in A, \; b\in B.$$
			Clearly, $\pi(A\otimes_\gamma B)\subseteq cl(\pi_1(A)\pi_2(B))$. Since $\pi$ is irreducible, 
			$$(\pi_1(A)\pi_2(B))'=cl(\pi_1(A)\pi_2(B))'\subseteq \pi(A\otimes_\gamma B)'=\mathbb{C}I.$$
			For the von Neumann algebra $\pi_1(A)''$, we have
			$$ Z(\pi_1(A)'') =\pi_1(A)''\cap \pi_1(A)'''=(\pi_1(A)'\cup \pi_1(A))'\subseteq (\pi_2(B)\cup \pi_1(A))' \subseteq \mathbb{C}I, $$	
				so that $\pi_1$ is a factor representation. Similarly, $\pi_2$ is a factor representation.
		\end{proof}
		\begin{thm}\label{Primitiveideal}
			For $C^\ast$-algebras $A$ and $B$, we have the following:
			\begin{enumerate}
				\item If $P \in Prim(A)$ and $Q \in Prim(B)$, then $A\otimes_\gamma Q + P\otimes_\gamma B \in Prim(A \otimes_\gamma B)$.
				\item If $K \in Prim(A \otimes_\gamma B)$, then $K= A\otimes_\gamma Q + P\otimes_\gamma B$ for some $P \in Prime(A)$ and $Q \in Prime(B)$.
				\item If $A$ and $B$ are seperable and $K \in Prim(A\otimes_\gamma B)$, then $K= A\otimes_\gamma Q + P\otimes_\gamma B$ for some $P \in Prim(A)$ and $Q \in Prim(B)$. 
			\end{enumerate}
		\end{thm}
		\begin{proof}
			(1) Since $P$ and $Q$ are primitive ideals, there exist irreducible *-representations $\pi_1$ and $\pi_2$ of $A$ and $B$ on Hilbert spaces $H_1$ and $H_2$, respectively,  such that $P=\ker \pi_1$ and $Q=\ker \pi_2$. Define $\pi : A\otimes B \rightarrow B(H_1 \otimes H_2)$ by				$\pi(a\otimes b)=\pi_1(a)\otimes \pi_2(b)$. Clearly $\pi$ is bounded with respect to $\otimes_{\gamma}$-norm, being  bounded with respect to $\omin$-norm, thus $\pi$ can be extended continuously  to $A\otimes_\gamma B$  as a *-representation. We first claim that $\pi$ is irreducible. It is easy to see that
			$$\pi(A\otimes_\gamma B)'\subseteq (\pi_1(A)\otimes \pi_2(B))'= (\pi_1(A)\overline\otimes \pi_2(B))'$$ where $\overline{\otimes}$ denotes the weak closure.  By Double
				Commutant Theorem, $\pi_1(A)$ and $\pi_2(B)$ are weakly dense in $\pi_1(A)''$ and $\pi_2(B)''$, so that                                                                                                                                                     $\pi_1(A)\overline{\otimes} \pi_2(B)=\pi_1(A)''\overline{\otimes} \pi_2(B)''$, thus $(\pi_1(A)\overline\otimes \pi_2(B))'=\pi_1(A)'\overline\otimes \pi_2(B)'\subseteq \mathbb{C}I$, by Tomita's Commutation Theorem, which shows that $\pi$ is irreducible. We claim that $\ker \pi=A\otimes_\gamma Q + P\otimes_\gamma B=K$ (say). Note that  $K=\ker q\subseteq \ker \pi$ where $q:A\otimes_\gamma B\rightarrow (A/P)\otimes_\gamma (B/Q)$ is the quotient map, by  \cite[Proposition 4]{GuptaJain Banach}. Therefore  $q(\ker \pi)$ is a closed ideal in $(A/P)\otimes_\gamma (B/Q)$, by \cite[Lemma 1]{GuptaJain Banach}.  If $(a+P)\otimes (b+Q) \in q(\ker \pi)$, then $a\otimes b\in \ker \pi$, which further gives either $\pi_1(a)=0$ or $\pi_2(b)=0$. Thus, either $a\in P$ or $b\in Q$, showing that $q(\ker \pi)$ has no non zero elementary tensor. Hence, by \cite[Corollary 2]{GuptaJain Banach}, $q(\ker \pi)= \{0\}$, so that $\ker \pi=K$.
				
				(2) Let $K=\ker \pi$, where $\pi$ is an irreducible $*$-representation of $A\otimes_\gamma B$. By \Cref{Factorrepresent}, $\pi (a\otimes b)=\pi_1(a)\pi_2(b)$ for $ a\in A$ and $b\in B,$ where $\pi_1$ and $\pi_2$ are factor representations of $A$ and $B$ on $H_1$ and $H_2$, respectively, with commuting ranges. 
				
				Let $P=\ker \pi_1$ and $Q=\ker \pi_2$. By \cite[II.6.1.11]{Blackadar}, $P\in Prime(A)$ and $Q \in Prime(B)$. It is easy to verify that $A\otimes_\gamma Q + P\otimes_\gamma B\subseteq K$. For any $a\otimes b\in K$, we have $\pi_1(a)\pi_2(b)=0$. Since $\pi_1(A)''$ is a factor and $\pi_2(B)''\subseteq \pi_1(A)'$, using \cite[Proposition IV.4.20]{Takesaki}, we deduce that either $\pi_1(a)=0$ or $\pi_2(b)=0$. It follows that $a\otimes b\in A\otimes_\gamma Q + P\otimes_\gamma B$. Then with the same argument as in first part we get $K=A\otimes_\gamma Q + P\otimes_\gamma B$.
				
				(3) If $A$ and $B$ are seperable $C^\ast$-algebras then every prime ideal is primitive. Hence, the proof follows from $(2)$.
		\end{proof}
		
	\section{Weak Centrality of Tensor Products}
	 In this section we discuss weak centrality of Haagerup and Banach space projective tensor product of $C^\ast$-algebras. We first introduce the concept of  weak centrality of Banach algebras.
	
	\begin{defn}
		A Banach algebra $X$ is said to be {\it weakly central} if the following two conditions are satisfied:
		\begin{enumerate}
			\item  No maximal modular ideal of $X$ contains $Z(X)$.
			\item For each pair of maximal modular ideals $M_1$ and $M_2$ of $X$, $M_1\cap Z(X)=M_2\cap Z(X)$ implies $M_1=M_2$.
		\end{enumerate}
			\end{defn}
	It is easy to observe that  condition $(1)$ is redundant, if $X$ is unital. 
			 Let us look at some examples of weakly and non-weakly central Banach algebras:

	\begin{example}
		Consider the Banach algebra $X:=\left\{\left(\begin{smallmatrix}
			S&T\\
			0&S
		\end{smallmatrix}\right) | S,T\in B(H) \right\}$ with $\| \left(\begin{smallmatrix}	S&T\\		0&S		\end{smallmatrix}\right)\|:=\|S\|_{B(H)}+\|T\|_{B(H)}$, where $H$ is an infinite dimensional seperable Hilbert space. Here, $Z(X)=\left\{\left(\begin{smallmatrix}			S&T\\			0&S		\end{smallmatrix}\right) | S,T\in \mathbb{C}I_{B(H)} \right\}$. Consider the closed ideal $M=\left\{\left(\begin{smallmatrix}
		S&T\\
		0&S
		\end{smallmatrix}\right) | S\in K(H),T\in B(H) \right\}$ of $X$. Since $K(H)$ is the only non-trivial proper closed ideal of $B(H)$, it is easy to verify that   $M$ is the only maximal modular ideal of $X$ with $\left(\begin{smallmatrix}
			I&0\\
			0&I
		\end{smallmatrix}\right)$ being the identity modulo $M$. Clearly $Z(X)\not\subseteq M$ and thus, $X$ is weakly central.
	\end{example}
	\begin{example}
		Consider $X:=\left\{\left(\begin{smallmatrix}
			a&b\\
			0&0
		\end{smallmatrix}\right)| a,b\in \mathbb{C} \right\}$ with any matrix norm. Here $Z(X)=\{0\}$ and $M=\left\{\left(\begin{smallmatrix}
			0&b\\
			0&0
		\end{smallmatrix}\right)| b\in \mathbb{C} \right\}$ is the only (maximal) modular ideal of $X$ with $\left(\begin{smallmatrix}
			1&0\\
			0&0
		\end{smallmatrix}\right)$ being the identity modulo $M$. Clearly $Z(X)\subseteq M$ and thus $X$ is not weakly central.
	\end{example}
	
	\begin{example}
		Take the Banach algebra $X:=\left\{\left(\begin{smallmatrix}
			a&0&d\\
			0&b&0\\
			0&0&c
		\end{smallmatrix}\right)| a,b,c,d\in \mathbb{C} \right\}$. It is easy to check that $Z(X)=\left\{\left(\begin{smallmatrix}
			t&0&0\\
			0&b&0\\
			0&0&t
		\end{smallmatrix}\right)| b,t\in \mathbb{C} \right\}$. Consider two distinct maximal modular ideals $$M_1=\left\{\left(\begin{smallmatrix}
			a&0&d\\
			0&b&0\\
			0&0&0
		\end{smallmatrix}\right)| a,b,d\in \mathbb{C} \right\}\ \text{and}\  M_2=\left\{\left(\begin{smallmatrix}
			0&0&d\\
			0&b&0\\
			0&0&c
		\end{smallmatrix}\right)| b,c,d\in \mathbb{C} \right\}$$ where $I$ (identity matrix) is the identity modulo for $M_1$ and $M_2$. It is easy to see that $M_1\cap Z(X)=M_2\cap Z(X)=\left\{\left(\begin{smallmatrix}
			0&0&0\\
			0&b&0\\
			0&0&0
		\end{smallmatrix}\right)| b\in \mathbb{C} \right\}$. Thus, $X$ is not weakly central. 
		\end{example}
	We next prove some interesting facts about weak centrality of Banach algebra. Before that we prove a basic result:
	\begin{lem}\label{maximalcentral}
		Let $X$ be a Banach *-algebra and $P$ be a primitive ideal of $X$ such that $Z(X)\not \subseteq P$. Then $P\cap Z(X)\in Max(Z(X))$.
	\end{lem}
	\begin{proof}
		Let $\pi$ be an irreducible $*$-representation of $X$ such that $\ker \pi=P$. Then, by \cite[Theorem 9.6.1]{Palmer}, $\pi(X)'=\mathbb{C} I$. Clearly, $\pi(Z(X))\subseteq \pi(X)'=\mathbb{C} I$. Since $Z(X)\not \subseteq P$, there exists $z\in Z(X)$ such that $z\not \in P$. Let, $\pi(z)=\alpha I$ for some $\alpha\neq 0$. For $\beta \in \mathbb{C}$,  $\beta z/\alpha\in Z(X)$ and $\beta I=\pi(\beta z/\alpha)$ so that $\mathbb{C} I\subseteq \pi (Z(X))$. Thus, $\pi|_{Z(X)}:Z(X)\rightarrow \mathbb{C} I$ is an onto map with $\ker \pi|_{Z(X)}=P\cap Z(X)$,  so that $Z(X)/(P\cap Z(X))\cong \mathbb{C} I$. Hence, $P\cap Z(X) \in Max(Z(X))$.  
	\end{proof}
	\begin{lem}\label{appid}
		Let $X$ be a Banach algebra and $I$ be an ideal of $X$ having an approximate identity. Then $Z(I)=I\cap Z(X)$.
	\end{lem}
	\begin{proof}
		Let $\{e_{\lambda}\}$ be an approximate identity of $I$ and $a\in Z(I)$. Then, for $x\in X$,  $xa=\lim x e_{\lambda} a=\lim a x e_{\lambda} = ax$, so that $a\in I \cap Z(X)$. Thus $Z(I)=Z(X)\cap I$.
	\end{proof}
		\begin{propn}
		Let $X$ be a non-unital Banach *-algebra and $\tilde{X}$ be its unitization. Then, $X$ is weakly central if and only if $\tilde{X}$ is weakly central.
	\end{propn}
	\begin{proof}
		Let $X$ be weakly central. If $\tilde{X}$ is not weakly central, 
		then there exist distinct ideals $M_1,M_2\in Max(\tilde{X})$ such that $M_1\cap Z(\tilde{X})=M_2\cap Z(\tilde{X})$. If $M_1$ and $M_2$ are both different from $X$, set $N_i=M_i\cap X$, $i=1,2$. By \cite[Chapter IV, \S 20C]{Loomis}, the map $Max_{X}(\tilde{X}) \ni N \to N\cap X \in Max(X)$ is a bijection, where $Max_{X}(\tilde{X}):=\{N \in Max(\tilde{X}): X \nsubseteq N\}$. Thus, $N_i\in Max(X)$ and clearly $N_1\cap Z(X)=N_2\cap Z(X)$ as $Z(X) \subseteq Z(\tilde{X})$. By weak centrality of $X$, $N_1=N_2$, and thus $M_1=M_2$, which is a contradiction. Now, suppose one of the ideals, say $M_1$, be equal to $X$. Since, $Z(\tilde{X})=Z(X)\oplus \mathbb{C}$, we have $Z(X)=X\cap Z(\tilde{X})=M_2\cap Z(\tilde{X})$.  Thus $Z(X)\subseteq M_2\cap X$, which contradicts the weak centrality of $X$ as $M_2\cap X \in Max(X)$. Hence, $\tilde{X}$ is weakly central.
		
	Conversely, let $\tilde{X}$ be weakly central.  Let $M\in Max(X)$ and if possible, let $Z(X) \subseteq M$. Now, $M=N\cap X$ for some $N\in Max_{X}(\tilde{X})$. Since $ X \in Max(\tilde{X})$, by weak centrality of $\tilde{X}$, $N\cap Z(\tilde{X})\neq X\cap Z(\tilde{X})=Z(X)$. Also, every maximal ideal is primitive, therefore by \Cref{maximalcentral}, we have $N\cap Z(\tilde{X})\in Max(Z(\tilde{X}))$. Now, $Z(X) \subset N \cap Z(\tilde{X}) \subset Z(\tilde{X})$, both the containments being proper, which is not true since $Z(X) \in Max(Z(\tilde{X}))$.  Thus,  $Z(X)\nsubseteq M$.

		Now suppose $M_1,M_2\in Max(X)$ such that $M_1\cap Z(X)=M_2\cap Z(X)$. Then there exist $N_i\in Max_X(\tilde{X})$ such that $M_i=N_i\cap X$, $i=1,2$. Since $Z(X)\not \subseteq M_i$, we have  $Z(X)\not \subseteq N_i$, so that $Z(X)\not \subseteq N_i\cap Z(\tilde{X})$. Again, by \Cref{maximalcentral}, $N_i\cap Z(\tilde{X}) \in Max(Z(\tilde{X}))$. Further
		\begin{eqnarray*}
			\notag
			(N_1 \cap Z(\tilde{X}))\cap Z(X)&=&(N_1 \cap X\cap Z(X))\cap Z(\tilde{X})\\
			\notag
			&=&(N_2 \cap X\cap Z(X))\cap Z(\tilde{X})\\
			&=&(N_2 \cap Z(\tilde{X}))\cap Z(X)	
		\end{eqnarray*}
		   Thus, $(N_1\cap Z(\tilde{X})) Z(X)\subseteq (N_2\cap Z(\tilde{X}))$. Since, $N_2\cap Z(\tilde{X})$ is a prime ideal of $Z(\tilde{X})$ being maximal and $Z(X)\not \subseteq N_2\cap Z(\tilde{X})$, we get $N_1\cap Z(\tilde{X})\subseteq N_2\cap Z(\tilde{X})$. Similarly, $N_2\cap Z(\tilde{X})\subseteq N_1\cap Z(\tilde{X})$. Since, $\tilde{X}$ is weakly central, this gives $N_1=N_2$ which would further imply that $M_1=M_2$. Hence the proof.
	\end{proof}
	
	The next result demonstrates that the quotient of a Banach algebra also respects the weak centrality. 
		\begin{propn}
		Every quotient of a weakly central Banach algebra is weakly central.
	\end{propn}
	\begin{proof}
		Let $X$ be weakly central and $I$ be any closed ideal of $X$. Consider $M/I\in Max(X/I)$ where $M\in Max(X)$ and $I\subseteq M$. Since $X$ is weakly central, $Z(X) \nsubseteq M$, which implies
		$(Z(X)+I)/I \nsubseteq M/I$. But $(Z(X)+I)/I \subseteq Z(X/I)$, thus $Z(X/I) \nsubseteq M/I$.
		
		Now, let $M_1/I,M_2/I\in Max(X/I)$ such that $M_1/I\neq M_2/I$. Since $X$ is weakly central and $M_1\neq M_2$ we have $M_1\cap Z(X)\neq M_2\cap Z(X)$, so there exists $x\in M_1\cap Z(X)$ such that $x\notin M_2\cap Z(X)$.  Then $x+I\in (M_1/I)\cap (Z(X)+I)/I$. However, $x+I\notin M_2/I$, 
		so that $(M_1/I)\cap Z(X/I)\neq (M_2/I)\cap Z(X/I)$. Hence, $X/I$ is weakly central.		 
	\end{proof}
	
		We are now ready to discuss our main results on the weak centrality for the Haagerup and Banach space projective tensor product of $C^\ast$-algebras. 	We first state two facts which will be needed for the further discussions. For the Banach space projective tensor product, the proofs of these results can be found in  \cite[Lemma 4, Theorem 7]{GuptaJain Banach}. For the Haagerup tensor product, the first part appears in \cite{ArchboldKaniuthSchlitingSomerset}. Through out the remaining article, $\oa$ represents either the Haagerup product or the Banach space projective product. 
		\begin{lem}\label{tensornormequality}
		Let $A$ and $B$ be $C^\ast$-algebras. 
		\begin{enumerate}
			\item  If $I_k$ and $J_k$, $k=1,2$ are $C^\ast$-subalgebras of $A$ and $B$ respectively such that $I_1\otimes_\alpha B+A\otimes_\alpha J_1= I_2\otimes_\alpha B+A\otimes_\alpha J_2$ then $I_1= I_2$ and $J_1= J_2$.
			\item $M \in Max(A \oa B) $ if and only if $M=M_1 \oa B + A \oa M_2$ for some $M_1 \in Max(A)$ and $M_2 \in Max(B)$.
		\end{enumerate}
				\end{lem}
		\begin{proof}
			 We only need to prove (2) for the Haagerup tensor product. Let $M\in Max(A\otimes_h B)$. Since $M$ is a maximal ideal of $A \oh B$, by \cite[Theorem 5.6]{AllenSinclairSmith}  $M=M_1 \otimes_h B + A \otimes_h M_2$ for some maximal ideals $M_1$ and $M_2$ of $A$ and $B$, respectively. As $(A\otimes_h B)/M\cong (A/M_1)\otimes_h (B/M_2)$ \cite[Corollary 2.6]{AllenSinclairSmith} and $(A\otimes_h B)/M$ is unital, thus $(A/M_1)\otimes_h (B/M_2)$ is also unital. This would further imply that  $A/M_1$ and $B/M_2$ are unital \cite[Theorem 1]{Loy}. Thus, $M_1$ and $M_2$ are modular.
			 
				Conversely, if $M=M_1 \otimes_h B + A \otimes_h M_2$, where $M_1\in Max(A)$ and $M_2\in Max(B)$, then by \cite[Theorem 5.6]{AllenSinclairSmith} $M$ is maximal in $A\otimes_h B$. Since, $A/M_1$ and $B/M_2$ are unital, using the above isomorphism, $(A\otimes_h B)/M$ is also unital and hence $M\in Max(A\otimes_h B)$.			
		\end{proof}

	\begin{propn}\label{maxidealcentre}
		Let $P=P_1\oa B + A\oa P_2\in Prim(A\oa B)$, where $P_1\in Prim(A)$ and $P_2\in Prim(B)$. Then, $$P \cap (Z(A \oa B))=(P_1\cap Z(A))\oa Z(B)+Z(A)\oa (P_2\cap Z(B)).$$
	\end{propn}
	\begin{proof}
		It is well known that  $Z(A\op B)=Z(A)\op Z(B)$ \cite[Theorem 2]{GuptaJain Banach} and  $Z(A\oh B)=Z(A)\oh Z(B)$ \cite[Theorem 2.13]{AllenSinclairSmith}. So, if $Z=Z(A)\oa Z(B)\subseteq P$ then the equality is trivially true. Now, let $Z\not\subseteq P$, then $Z(A)\not\subseteq P_1$ and $Z(B)\not\subseteq P_2$. By \Cref{maximalcentral}, $P_1\cap Z(A)\in Max(Z(A))$ and $P_2\cap Z(B)\in Max(Z(B))$. Since Haagerup and Banach space projective tensor product both are injective injective for $C^\ast$-subalgebras \cite[Theorem 4.4]{PaulsenSmith}, \cite[Theorem 1]{GuptaJain Banach}, we have
		$$(P_1\cap Z(A))\otimes_\alpha Z(B)\subseteq(P_1\otimes_\alpha B)\cap (Z(A)\otimes_\alpha Z(B))\subseteq P\cap Z,$$
		$$Z(A)\otimes_\alpha(P_2\cap Z(B))\subseteq (Z(A)\otimes_\alpha Z(B)) \cap (A \oa P_2) \subseteq P\cap Z.$$
		As a consequence, 
		$$(P_1\cap Z(A))\otimes_\alpha Z(B)+Z(A)\otimes_\alpha(P_2\cap Z(B))\subseteq P\cap Z.$$
		Thus, by \Cref{tensornormequality}, $(P_1\cap Z(A))\oa Z(B)+Z(A)\oa(P_2\cap Z(B))$ is a maximal ideal of $Z$. Since, $Z\not\subseteq P$, $P\cap Z$ is a proper ideal of $Z$. Therefore, $$(P_1\cap Z(A))\otimes_\alpha Z(B)+Z(A)\otimes_\alpha(P_2\cap Z(B))= P\cap Z,$$ and hence the result.
	\end{proof}
With all the ingredients prepared, we now prove our main result.
	\begin{thm}\label{weaklycentral}
			If $A$ and $B$ are weakly central, then so is $A\otimes_\alpha B$. The converse is true if $Max(A)$ and $Max(B)$ both are non-empty.
	\end{thm}
	\begin{proof}
		Consider a maximal modular ideal $M$ of $A \oa B$, then by \Cref{tensornormequality}, $M=M_1\otimes_\alpha B+A\otimes_\alpha M_2$ for some $M_1\in Max(A)$ and $M_2\in Max(B)$. Since $A$ and $B$ are weakly central, $Z(A)\not\subseteq M_1$ and $Z(B)\not\subseteq M_2$. So, there exist $z_1\in Z(A)$ and $z_2\in Z(B)$ such that $z_1\notin M_1$ and $z_2\notin M_2$. Choose $\phi \in A^*$ and $\psi \in B^*$ such that $\phi (z_1)\neq 0$, $\psi (z_2)\neq 0$ and $\phi(M_1)=\psi(M_2)=0$. Then, $\phi\otimes_\alpha \psi \in (A\otimes_\alpha B)^*$ with $(\phi\otimes_\alpha \psi) (M)=0$ and $(\phi\otimes_\alpha \psi)(z_1\otimes z_2)\neq 0$, so that $z_1\otimes z_2\in Z(A \oa B)$ but $z_1\otimes z_2 \notin M$, giving that $Z(A\otimes_\alpha B)\nsubseteq M$. Thus, the first condition is satisfied. 
		
		Next, consider $M,M'\in Max(A\otimes_\alpha B)$ such that $M\cap Z=M'\cap Z$, where $Z=Z(A \oa B)= Z(A)\otimes_\alpha Z(B)$. Now, $M=M_1\otimes_\alpha B+A\otimes_\alpha M_2$ and $M'=M'_1\otimes_\alpha B+A\otimes_\alpha M'_2$ for some $M_1,M_1'\in Max(A)$ and $M_2,M_2'\in Max(B)$. Also, as done above, $Z \nsubseteq M$ and $Z \nsubseteq M'$, so $M\cap Z$ and $M'\cap Z$ are proper ideals of $Z$. Therefore, by \Cref{maxidealcentre},
		$$M\cap Z=(M_1\cap Z(A))\otimes_\alpha Z(B)+Z(A)\otimes_\alpha(M_2\cap Z(B)),$$
		$$M'\cap Z=(M'_1\cap Z(A))\otimes_\alpha Z(B)+Z(A)\otimes_\alpha(M'_2\cap Z(B)).$$
		Since $M\cap Z=M'\cap Z$, by \Cref{tensornormequality}, $M_1\cap Z(A)=M'_1\cap Z(A)$ and $M_2\cap Z(B)=M'_2\cap Z(B)$. Since $A$ and $B$ are weakly central, we have $M_i=M'_i$, for $i=1,2$. Hence $M=M'$ and $A\oa B$ is weakly central.
		
		For the converse, let if possible, $A$ be not weakly central. Then there are two possibilities:
		
		Case (1): There exists $M_1\in Max(A)$ such that $Z(A)\subseteq M_1$. Consider $M=M_1\otimes_\alpha B+A\otimes_\alpha M_2$ for some $M_2\in Max(B)$. Then, by \Cref{tensornormequality}, $M\in Max(A\otimes_\alpha B)$. Thus, $Z(A\otimes_\alpha B) \subseteq M$, which contradicts the fact that $A\otimes_\alpha B$ is weakly central. 
		
		Case (2): There exist $M_1,M'_1\in Max(A)$,  $M_1 \neq M'_1$, such that $M_1\cap Z(A)=M'_1\cap Z(A)$. For some fixed $M_2\in Max(B)$ take $M=M_1\otimes_\alpha B+A\otimes_\alpha M_2$ and $M'=M'_1\otimes_\alpha B+A\otimes_\alpha M_2$. Then, by \Cref{tensornormequality} $M,M'\in Max(A\otimes_\alpha B)$ and $M\neq M'$. Further, $M\cap Z=(M_1\cap Z(A))\otimes_\alpha Z(B)+Z(A)\otimes_\alpha(M_2\cap Z(B))$ and $M'\cap Z=(M'_1\cap Z(A))\otimes_\alpha Z(B)+Z(A)\otimes_\alpha(M_2\cap Z(B))$. This gives that $M\cap Z=M'\cap Z$, which is a contradiction to the weakly centrality of $A\otimes_\alpha B$. Hence, $A$ is weakly central, and similarly, $B$ is weakly central.
	\end{proof}
	\begin{remark}\label{twc}
		The converse of \Cref{weaklycentral} is not true if we drop the conditions that $Max(A)$ and $Max(B)$ are non empty. In fact, if any one of these spaces is empty, then by \Cref{tensornormequality}, $Max(A \oa B) = \emptyset$, thus $A\oa B$ is always weakly central irrespective of the weak centrality of $A$ or $B$.
	\end{remark}
	In a unital algebra, every ideal is modular. Thus, in the case of unital $C^*$-algebra, we have the following characterization.
	\begin{cor}
		If $A$ and $B$ are unital $C^\ast$-algebras, then $A\oa B$ is weakly central if and only if $A$ and $B$ are weakly central.
	\end{cor}

		
	
	We next investigate the largest weakly central ideal of $A \oa B$, denoted by $J_{wc}(A\otimes_\alpha B)$ (if exists). For this purpose, let us first fix some notations. For a Banach algebra $X$, we set
	\begin{itemize}
		\item $T_X^1:=\{M \in Max(X)$: $Z(X)\subseteq M$\}.
\item  	$T_X^2:=\{ M\in Max(X) \setminus T_X^1 : M\cap Z(X)=N\cap Z(X) \ \text{for some} \  N\in Max(X) \setminus \{M\}\}$.
		\item $T_X=T_X^1\cup T_X^2$.
	\end{itemize}
	Clearly, $X$ is weakly central if and only if $T_X = \phi$. Thus, $T_X$ contains those maximal modular ideals which prevent $X$ to be weakly central.
		
	\begin{propn}\label{Banachlargestweaklycentral}
		Let $X$ be a Banach algebra and $I$ be a weakly central ideal of $X$ possessing an approximate identity. Then $I \subseteq \ker T_X$.
	\end{propn}
	
	
	\begin{proof} 
 Let, if possible, $I\not \subseteq \ker T_X$, then either $I\not \subseteq \ker T_{X}^{1}$ or $I\not \subseteq \ker T_{X}^{2}$. If $I\not \subseteq \ker T_{X}^{1}$, then there exists $M\in T_{X}^{1}$ such that $I\not \subseteq M$. By \cite[Proposition 1.1]{Barnes} $M\cap I\in Max(I)$. Since $M\in T_{X}^{1}$, we have $Z(X)\subseteq M$. By \Cref{appid} $Z(I)=Z(X)\cap I\subseteq M\cap I$ which contradicts the weak centrality of $I$. 
		
		If $I\not \subseteq \ker T_{X}^{2}$, then there exists $M_1\in T_{X}^{2}$ such that $I\not \subseteq M_1$. Since $M_1 \in T_{X}^{2}$ there exists $M_2\in Max(X)$ such that $M_1 \neq M_2$, $Z(X) \nsubseteq M_2$ and $M_1\cap Z(X)=M_2\cap Z(X)$. Then, 
		\begin{eqnarray}\label{2}
			\notag
			(M_1\cap I)\cap Z(I)&=&M_1\cap I\cap Z(X) \cap I\\
			\notag
			&=&M_2\cap Z(X)\cap I\\
			&=&(M_2\cap I)\cap Z(I).
		\end{eqnarray}
		This implies that $I\not \subseteq M_2$, otherwise by \eqref{2}, $Z(I)\subseteq M_1\cap I$ which contradicts the weak centrality of $I$ since $M_1\cap I\in Max(I)$. Thus, $ M_2\cap I\in Max(I)$ and weak centrality of $I$ together with \eqref{2} imply that $M_1\cap I=M_2\cap I$. Since, $M_1\cap I\subseteq M_2$, $M_2$ is a prime ideal and $I\not \subseteq M_2$,we have $M_1\subseteq M_2$. Similarly, $M_2\subseteq M_1$, so that $M_1=M_2$ which is a contradiction and hence the result.  
	\end{proof}
	
	\begin{cor}\label{largestweaklycentralhaagerup}
		Let $A$ and $B$ be $C^\ast$-algebras with either of them having finitely many closed ideals. Then every weakly central ideal of $A\oa B$ is contained in $\ker T_{A\oa B}$.
	\end{cor}
	\begin{proof}
		Since one of $A$ or $B$ has finitely many closed ideals, every closed ideal of $A\oa B$ is a sum of finitely many product ideals, by \cite[Theorem 5.3]{AllenSinclairSmith} and \cite[Theorem 5]{GuptaJain Banach}. Thus, by \cite[Lemma 3.3]{AllenSinclairSmith} and \cite[Corollary 5]{GuptaJain Banach}, every closed ideal of $A\oa B$ has an approximate identity. The result now follows from \Cref{Banachlargestweaklycentral}.
	\end{proof}
	
	For a $C$*-algebra $A$, Archbold and Gogi\'{c}  proved that $\ker T_A$ is the largest weakly central ideal of $A$ and thus coincides with $J_{wc}(A)$  \cite[Theorem 3.22]{ArchboldGogic}. We investigate the cases for which the same is true for the tensor products. In order to do that, we first derive an explicit expression for $\ker T_{A\oa B}$ in terms of $\ker T_A$ and $\ker T_B$.
	
	\begin{lem}\label{largestweaklycentral}
		Let $A$ and $B$ be $C$*-algebras. Then
		$$\ker T_{A\oa B}=\mathcal{M}(A)\oa B+A\oa\mathcal{M}(B) + \ker T_A\oa \ker T_B.$$
	\end{lem}
	\begin{proof}
		Let us write $X = A \oa B$. We first claim that 
		$T_{X}^1=S_1\cup S_2$ , where \[S_1=\{M\oa B + A\oa M'\;| \; M\in T_A^1, M'\in Max(B)\}\]
		\[S_2=\{M\oa B + A\oa M'\ | M\in Max(A), M'\in T_B^1\}\]
		
		 Consider $I=M\oa B + A\oa M'\in S_1\cup S_2 $ where either $M\in T_A^1$ or $M'\in T_B^1$. Then, either  $Z(A)\subseteq M$ or  $Z(B)\subseteq M'$ which gives $Z(A)\oa Z(B)\subseteq I$. Since $I \in Max(X)$, we have $ I \in T_{X}^1$.	
		For the reverse inclusion, if $I \in T_{X}^1$, then $I=M\oa B + A\oa M'$, for some $M\in Max(A)$ and $M'\in Max(B)$ with $ Z(A)\oa Z(B)\subseteq I$. 
	If $Z(A)\nsubseteq M$ and $Z(B)\nsubseteq M'$ then for $z_1\in Z(A) \setminus M$ and $z_2\in Z(B) \setminus M'$, choose $\phi \in A^*$ and $\psi \in B^*$ such that $\phi(z_1)\neq 0$ and $\psi(z_2)\neq 0$ but $\phi(M)=\psi(M')=0$. Then $\phi \otimes \psi \in (A\oa B)^*$ with $(\phi \otimes \psi)(z_1\otimes z_2)\neq 0$ and $(\phi \otimes \psi)(I)=0$. Thus $z_1\otimes z_2\not\in I$ which is a contradiction. Hence, either $Z(A)\subseteq M$ or $Z(B)\subseteq M'$, and we have the equality.
	
		Next, we claim that 
		 $T_{X}^2=S_3\cup S_4$, where 
		$$ S_3=\{M\oa B + A\oa M'|\ M\in T_A^2,\  M' \in Max(B) \setminus T_B^1\}$$
		$$S_4=\{M\oa B + A\oa M'|\  M\in Max(A) \setminus T_A^1, \ M'\in T_B^2\}$$
		Let $I=M\oa B + A\oa M'\in S_3$, so that $M\not\in T_A^1$ and $M' \not\in T_B^1$. Then, as done above, $Z(A)\oa Z(B)\not\subseteq I$. Since $M\in T_A^2$, there exists $N\in Max(A)$ such that $M\neq N$ and $M\cap Z(A)=N\cap Z(A)$. Set $J=N\oa B+A\oa M'$, then by \Cref{tensornormequality}, $I\neq J$. Also, by \Cref{maxidealcentre} and the fact that every maximal ideal is a primitive ideal, we have
		 \begin{eqnarray*}
		 	I\cap (Z(A)\oa Z(B))&=&(M\cap Z(A))\oa Z(B)+Z(A)\oa (M'\cap Z(B))\\
		 	&=&(N\cap Z(A))\oa Z(B)+Z(A)\oa (M'\cap Z(B))\\
		 	&=&J\cap (Z(A)\oa Z(B)).
		 \end{eqnarray*}
		Thus, $I\in T_{X}^2$. The other case follows on the similar lines. 
		 
		For the reverse inclusion, let $I\in T_{X}^2$, then $I=M\oa B + A\oa M'$ for some $M\in Max(A)$ and $M'\in Max(B)$. Since, $Z(A)\oa Z(B)\not\subseteq I$, therefore $M\not\in T_A^1$ and $M'\not\in T_B^1$. Also, there exists $J=N\oa B+A\oa N'\in Max(A\oa B)$ such that $I\neq J$ and
\begin{equation}\label{4}
	I\cap (Z(A)\oa Z(B))=J\cap (Z(A)\oa Z(B)).
\end{equation}  
Clearly, either $M\neq N$ or $M'\neq N'$. Also, using \Cref{maxidealcentre} in  \Cref{4}, we get 
\begin{eqnarray*}
		(M\cap Z(A))\oa Z(B)+Z(A)\oa (M'\cap Z(B))&=&(N\cap Z(A))\oa Z(B)\\
		& &+Z(A)\oa (N'\cap Z(B)).
\end{eqnarray*}

An application of \Cref{tensornormequality} now gives that either $M\in T_A^2$ or $M'\in T_B^2$ and hence the claim.
		Now, using \cite[Lemma 1.3]{ArchboldKaniuthSchlitingSomerset} and \Cref{kerideals}, we get
		 \begin{align*}		 	
		 Y_1 &:=  \ker S_1=\ker T_A^1 \oa B+A\oa \mathcal{M}(B), \\
		 Y_2 &:=  \ker S_2=\mathcal{M}(A)\oa B+A\oa \ker T_B^1, \\
		 Y_3 &:= \ker S_3=\ker T_A^2 \oa B+A\oa \ker(Max(B) \setminus T_B^1),\\
		Y_4 &:= \ker S_4=\ker(Max(A) \setminus T_A^1) \oa B+A\oa \ker T_B^2.
	\end{align*}
		Further, using distributivity of intersection of ideals as given in \cite[Proposition 3.6]{AllenSinclairSmith}, \cite[Corollary 4.6]{Smith} and \Cref{intersection}, we have
		$$Y_1\cap Y_3=\ker T_A\oa B+\ker T_A^1\oa \ker(Max(B) \setminus T_B^1)+A\oa \mathcal{M}(B),$$
		$$Y_2\cap Y_4=\mathcal{M}(A)\oa B+\ker(Max(A) \setminus T_A^1)\oa \ker T_B^1+A\oa \ker T_B.$$
		Finally,
		$$\ker T_{X}=\underset{i=1}{\overset{4}{\bigcap}} Y_i=\mathcal{M}(A)\oa B+A\oa \mathcal{M}(B)+\ker T_A\oa \ker T_B.$$
	\end{proof}
	
	 We next prove another main result of this article which assures that the largest weakly central ideal of $A\oa B$ coincides with $\ker T_{A\oa B}$ in certain cases. Recall that a $C^\ast$-algebra $A$ is said to be {\it strongly semisimple} if $\mathcal{M}(A)=\{0\}$.

	\begin{thm}\label{speciallwc}
		Let $A$ and $B$ be $C$*-algebras such that one of them contains finitely many closed ideals, and $Max(A\oa B)$ be non empty. Let any one the following conditions hold:
		\begin{enumerate}
			\item Either $A$ or $B$ is simple,
			\item Both $A$ and $B$ are strongly semisimple, 
			\item  Either $J_{wc}(A)=\{0\}$ or $J_{wc}(B)=\{0\}$.
		\end{enumerate}
		Then 
		$$J_{wc}(A\oa B)=\mathcal{M}(A)\oa B+A\oa\mathcal{M}(B) + J_{wc}(A)\oa J_{wc}(B)$$ 
		
	\end{thm}
	\begin{proof}
		By \Cref{largestweaklycentralhaagerup}, we know that every weakly central ideal of $A\oa B$ is contained in $\ker T_{A\oa B}$. Thus, it is sufficient to prove that  $\ker T_{A\oa B}$ is weakly central and the result will then follow from \Cref{largestweaklycentral}. 
		\begin{enumerate}
			\item Assume $B$ is simple, so that $\mathcal{M}(B)=\{0\}$ and $J_{wc}(B)=B$. Then
			$$\ker T_{A\oa B}=\mathcal{M}(A)\oa B+J_{wc}(A)\oa B=J_{wc}(A)\oa B,$$ 
			Thus, using \Cref{weaklycentral},  $\ker T_{A\oa B}$ is weakly central. 
			\item Since $A$ and $B$ are strongly semisimple, by \Cref{largestweaklycentral}
			 $$\ker T_{A\oa B}=J_{wc}(A)\oa J_{wc}(B),$$ which is weakly central by \Cref{weaklycentral}. 
		\item If  $J_{wc}(A)=\{0\}$, then $\mathcal{M}(A)=\{0\}$, so that
			$$\ker T_{A\oa B}=A\oa \mathcal{M}(B).$$
			By \cite[Lemma 2.2]{ArchboldGogic}, $\mathcal{M}(B)$ does not exhibit any proper modular ideal. Thus, by \Cref{twc}, $A\oa \mathcal{M}(B)$ is weakly central. \qedhere					
		\end{enumerate}
	\end{proof}
			We next illustrate this result with the help of an example.
	\begin{example}\label{eg1}
		Let $A$ be the rotation algebra or the full $C^\ast$-algebra of the free group $\mathbb{F}_2$ on two generators, then $J_{wc}(A)=\{0\}$ (see, \cite[Example 3.24, 3.25]{ArchboldGogic}). For $B=B(H)$, $H$ being infinite dimensional seperable Hilbert space, by \Cref{largestweaklycentral}, $\ker T_{A\oa B}=A\oa K(H)$, which is clearly weakly central since $ K(H)$ does not have any maximal modular ideal.
	\end{example} 
	Unlike in the case of $C^*$-algebras, it is not known if \Cref{speciallwc} holds in general for arbitrary $C^*$-algebras, however, we  provide one more instance in which $\ker T_{A\oh B}$ turns out to be the largest weakly central ideal of $A\oh B$.
	\begin{example}\label{eg2}
		Let $A$ be the $C^\ast$-algebra consisting of all continuous functions $f:[0,1]\rightarrow \mathbb{M}_2(\mathbb{C})$ such that $f(1)=diag(\lambda(f),0)$, for some scalar $\lambda(f)\in \mathbb{C}$.  		
		For $t\in [0,1]$, define $M_t=\{f\in A: f(t)=0\}$. By \cite[V.26.2]{Naimark}, $Max(A)=\{M_t: t\in [0,1]\}$. Note that $A$ is not weakly central since 
		$$Z(A)=\{diag(g,g): g\in C[0,1],g(1)=0\} \subseteq M_1.$$
		 Let  $B= B(H)$, where $H$ is an infinite dimensional seperable Hilbert space. We claim that $\ker T_{A\oh B}$ is weakly central, so that by \Cref{largestweaklycentralhaagerup}, $J_{wc}(A \oh B) =\ker T_{A\oh B}$.  Observe that $M_1$ is isomorphic to $C_0([0,1),\mathbb{M}_2(\mathbb{C}))$, the $C^*$-algebra consisting of all continuous functions $f:[0,1) \to \mathbb{M}_2(\mathbb{C})$ such that $\lim\limits_{t \to 1}f(t) =0 $. Thus, using \cite[Theorem 3.5]{TalwarJain} and the fact that $\mathbb{M}_2(\mathbb{C})$ is weakly central, we get $M_1$ is weakly central. In fact, the maximality of $M_1$ implies that it is the  largest weakly central ideal of $A$ since it is contained in $\ker T_A$, the largest weakly central ideal of $A$. By \Cref{largestweaklycentral}, and the fact that $\mathcal{M}(A)=\underset{t\in [0,1]}{\bigcap} M_t=\{0\}$, we have
		 \begin{eqnarray}\label{kerrel}
		 	\ker T_{A\otimes_h B}=A\otimes_h K(H)+ M_1\otimes_h B.
		 \end{eqnarray}
		We shall first determine the general form of any maximal modular ideal of $\ker T_{A\oh B}$.   Since  $\ker T_{A\oh B}$ possesses a bounded approximation identity, every closed ideal of $\ker T_{A\oh B}$ is an ideal of $A\oh B$, and thus  is a finite sum of product ideals, since $B$ has finitely many closed ideals. 
		Therefore, every proper closed ideal of $\ker T_{A\oh B}$ can be written as $P\oh K(H)+Q\oh B$ where $P$ is a closed ideal of $A$ and $Q$ is a closed ideal of $M_1$. If $M$ is a maximal closed ideal of $\ker T_{A\oh B}$, then either $M=M_1\oh B$ or $M=A\oh K(H)+ J\oh B$ where $J$ is maximal in $M_1$.  However, $M_1 \oh B$ is a non modular ideal as by Second Isomorphism Theorem, $$\frac{A\otimes_h K(H)+ M_1\otimes_h B}{M_1\otimes_h B}\cong \frac{A\otimes_h K(H)}{M_1\otimes_h K(H)}\cong (A/M_1)\otimes_h K(H),$$ which is non unital since $K(H)$ is non-unital. Now, using Second and Third Isomorphism Theorems, we have   \begin{eqnarray*}
			\frac{A\otimes_h K(H)+ M_1\otimes_h B}{A\oh K(H)+J\otimes_h B}&\cong& \frac{\left(\frac{A\otimes_h K(H)+ M_1\otimes_h B}{A\oh K(H)}\right)}{\left(\frac{A\oh K(H)+J\otimes_h B}{A\otimes_h K(H)}\right)}\\
			&\cong& \frac{\left(\frac{M_1\otimes_h B}{M_1\otimes_h K(H)}\right)}{\left(\frac{J\otimes_h B}{J\otimes_h K(H)}\right)}\\
			&\cong& \frac{\left(M_1\oh \frac{B}{K(H)}\right)}{\left(J\oh \frac{B}{K(H)}\right)}\\
			&\cong& \left(\frac{M_1}{J}\right)\oh \left(\frac{B}{K(H)}\right),
		\end{eqnarray*}
		which is unital if and only if $J$ is modular. Thus, $A\oh K(H)+J\otimes_h B$ is a maximal modular ideal of $\ker T_{A\oh B}$ if and only if $J$ is maximal modular in $M_1$. Further, by \cite[Lemma 2.2]{ArchboldGogic}, any maximal modular ideal of $M_1$ is of the form $M_t\cap M_1$ where $t\in [0,1)$. Hence, any maximal modular ideal of $\ker T_{A\oh B}$ is of the form $M^t=A\oh K(H)+ (M_t\cap M_1)\oh B$. Observe that, for $t \in [0,1)$, $M_t\cap M_1=\{f\in A: \ f(1)=f(t)=0\}$.  
		
		Next, we claim that $Z(\ker T_{A\oh B})=Z(A)\oh \C I$. Since $\ker T_{A\oh B}$ has an approximate identity, by \Cref{appid} and \cite[Corollary 4.6]{Smith}, we have $$Z(\ker T_{A\oh B}) = Z(A\oh B)\cap \ker T_{A\oh B} =(Z(A)\oh \C I) \cap \ker T_{A\oh B}.$$ 
		This, together with Equation \ref{kerrel}, gives that
		$$(Z(A)\oh \C I)\cap (A\otimes_h K(H))+ (Z(A)\oh \C I)\cap (M_1\otimes_h B) \subseteq Z(\ker T_{A\oh B}),$$  
		so that $Z(A)\oh \C I \subseteq Z(\ker T_{A\oh B})$ and hence the claim.
		
		Note that $Z(A)\not\subseteq M_t\cap M_1$  as for $f(x)=diag(x-1,x-1)$, $f\in Z(A)$ but $f\not\in M_t\cap M_1$. Thus, $Z(\ker T_{A\oh B})=Z(A)\oh \C I\not\subseteq M^t$. 
		
		Finally, consider  $ Z(\ker T_{A\oh B}) \cap M^s = Z(\ker T_{A\oh B}) \cap M^t$ for some $s\neq t, s,t \in [0,1)$. As done earlier, 
		 \begin{eqnarray*}
		Z(\ker T_{A\oh B}) \cap M^s & = & (Z(A) \oh \C I)\cap \big(A\oh K(H)+ (M_s\cap M_1)\oh B\big) \\
		 & = & (Z(A)\cap M_s)\oh \C I.
		 \end{eqnarray*}
	Thus, we have
		 $$ (Z(A)\cap M_s)\oh \C I =(Z(A)\cap M_t)\oh \C I,$$
which further implies $Z(A)\cap M_s = Z(A)\cap M_t$. Observe that $Z(A)\cap M_s=\{diag(g,g): g\in C[0,1],\ g(1)=g(s)=0\}$. If $M_s\neq M_t$, then for $f(x)=diag((x-s)(x-1)$, $f \in Z(A)\cap M_s$ but $f\not\in Z(A)\cap M_t$, so that $M_s = M_t$. Hence, $M^s = M^t$, and $\ker T_{A\oh B}$ is weakly central.
		
	\end{example}
	
	
	\section{Centrality and Quasi-centrality of tensor products}
	We first define the notion of centrality for Banach algebras in  general.
	\begin{defn}
		A Banach algebra $X$ is called {\it central} if the following two conditions are satisfied:
		\begin{enumerate}
			\item  No primitive ideal of X contains $Z(X)$.
			\item For each pair of primitive ideals $P$ and $P'$ of $X$, $P\cap Z(X)=P'\cap Z(X)$ implies $P=P'$.
		\end{enumerate}
	\end{defn}
	We characterize the centrality of $A\oa B$ in terms of centrality of $A$ and $B$.
		\begin{thm}\label{conversecentral}
		If $A\oa B$ is central then so are $A$ and $B$. The converse is true, if $A$ and $B$ both are separable. 
		\end{thm}
		\begin{proof}
		We first claim that $A$ is central. Let $P_1\in Prim(A)$ be such that $Z(A)\subseteq P_1$. For a fixed $P_2\in Prim(B)$, by \Cref{Primitiveideal} and \cite[Theorem 5.13]{AllenSinclairSmith}, $P=P_1\oa B+ A\oa P_2\in Prim(A\oa B)$, and $Z(A)\oa Z(B)\subseteq P$, which contradicts centrality of $A\oa B$. Thus, $Z(A)\not\subseteq P_1$. Next, consider $P_1,P_1'\in Prim(A)$ such that $P_1\cap Z(A)=P_1'\cap Z(A)$. For $P_2 \in Prim(B)$ set $P=P_1\oa B+ A\oa P_2$ and $P'=P_1'\oa B+ A\oa P_2$, then $P,P'\in Prim(A\oa B)$. By \Cref{maxidealcentre}, we have $P\cap Z(A\oa B)=P'\cap Z(A\oa B)$ which further gives $P = P'$, since $A\oa B$ is central. Finally, using \Cref{tensornormequality}, we get $P_1 = P_1'$ and hence $A$ is central. Similarly, $B$ is also central.
		
		For the converse, let $P$ be a primitive ideal  of $A\oa B$.	Since $A$ and $B$ are separable, by \Cref{Primitiveideal} and \cite[Theorem 5.13]{AllenSinclairSmith}, $P=P_1\oa B+A\oa P_2$ where $P_1\in Prim(A)$ and $P_2\in Prim(B)$.
		Clearly, $Z(A)\oa Z(B)\not\subseteq P$, otherwise, either $Z(A)\subseteq P_1$ or $Z(B)\subseteq P_2$ which is a contradiction. Consider, another primitive ideal $P'=P'_1\oa B+A\oa P'_2$   of $A\oa B$ such that $P\cap Z(A\oa B)=P'\cap Z(A\oa B)$, where $P_1' \in Prim(A)$ and $P_2' \in Prim(B)$. By \Cref{maxidealcentre},
		\begin{eqnarray*}
			(P_1\cap Z(A))\oa Z(B)+Z(A)\oa (P_2\cap Z(B))&=&(P'_1\cap Z(A))\oa Z(B)\\
			& &+Z(A)\oa (P'_2\cap Z(B)),
		\end{eqnarray*}		
	  which further imply $P_1\cap Z(A)=P'_1\cap Z(A)$ and $P_2\cap Z(B)=P'_2\cap Z(B)$, by \Cref{maximalcentral} and \Cref{tensornormequality}. Since $A$ and $B$ are central, we have $P_i=P'_i$ for $i=1,2$, and hence $P=P'$. Thus, $A\oa B$ is central.		
		\end{proof}
	
	We finally establish that quasi-centrality is passed to the tensor products. Recall that a Banach algebra is said to be {\it quasi-central} if no primitive ideal contains its center. Archbold (\cite[Proposition 1]{Archboldcentre}) proved that a $C^*$-algebra is quasi central if and only if it has a central approximation identity (approximate identity with each element in the centre). This is indeed true for the tensor products too. 
	\begin{thm}
		For $C$*-algebras $A$ and $B$, the following are equivalent:
		\begin{enumerate}
			\item $A\oa B$ is quasi-central.
			\item $A$ and $B$ are quasi-central.
			\item $A\oa B$ contains central approximate identity.
			\item Closed ideal generated by $Z(A)\oa Z(B)$ is $A\oa B$.
		\end{enumerate}
	\end{thm}
	\begin{proof}
		
		$(1)\implies (2)$ follows on the similar lines of \Cref{conversecentral}.\\
		$(2)\implies (3)$ Since $A$ and $B$ are quasi-central, by \cite[Proposition 1]{Archboldcentre}, $A$ and $B$ contain central approximate identities, say $\{e_\lambda\}_{\lambda \in C}$ and $\{f_\mu\}_{\mu \in D}$, respectively. Now consider the  set $C\times D$ with partial order as $(\lambda_1,\mu_1)\leq (\lambda_2,\mu_2)$ if and only if $\lambda_1\leq \lambda_2$ and $\mu_1 \leq \mu_2$. If we set $e_{(\lambda,\mu)}=e_\lambda$ and $f_{(\lambda,\mu)}=f_\mu$, then  $\{e_{(\lambda,\mu)} \otimes f_{(\lambda,\mu)}\}_{(\lambda,\mu)\in C\times D}\in Z(A)\oa Z(B)$ is a central approximate identity of $A\oa B$.\\
		$(3)\implies (4)$ This is trivial.\\
		$(4)\implies (1)$ Since the closed ideal generated by $Z(A)\oa Z(B)$ is $A\oa B$,  no proper primitive ideal contains $Z(A)\oa Z(B)$. Hence, $A\oa B$ is quasi-central.
	\end{proof}

\end{document}